\documentclass[reqno,12pt]{amsart}

\theoremstyle{remark}

\numberwithin{equation}{section}

\allowdisplaybreaks

\begin{document}

\title[Using  integrals of squares of  real-valued  functions ]
{Using integrals of squares of  certain real-valued special functions
 to prove that  the P\'olya  $\Xi^*(z)$ function, the  functions
$K_{iz}(a), a>0,$ and  some other
 entire functions have only real zeros}

\author{George Gasper}
\address{Department of Mathematics, Northwestern University,
2033 Sheridan Road, Evanston, IL 60208-2730, USA}
\email{george@math.northwestern.edu}
\urladdr{http://www.math.northwestern.edu/{\textasciitilde}george}

\subjclass\ {}

\dedicatory{Dedicated to Dan Waterman
on the occasion of his 80th birthday}

\begin{abstract}
Analogous to the use of sums of squares of  certain real-valued special functions 
to prove the reality of the zeros of the 
Bessel functions $J_\alpha (z)$
when $\alpha \ge -1,$ confluent hypergeometric functions ${}_0F_1(c\/; z)$
when $c>0$ or $0>c>-1$, Laguerre polynomials $L_n^\alpha(z)$ when 
$\alpha \ge -2,$ Jacobi polynomials $P_n^{(\alpha,\beta)}(z)$
when $\alpha \ge -1$ and $ \beta \ge -1,$ and some other entire special functions
 considered in
 G. Gasper [{\it Using sums of squares to prove that certain entire functions have 
only real zeros,\/}   in {\it Fourier 
Analysis: Analytic and Geometric Aspects, \/} W. O. Bray, P. S. Milojevi\'c and 
C. V. Stanojevi\'c, eds., Marcel Dekker, Inc., 1994, 171--186.],
 integrals of squares of  certain real-valued special functions  are used 
to prove the reality of the zeros of 
the  P\'olya $\Xi^*(z)$ function, the  $K_{iz}(a)$ functions
when $a>0,$
and some other entire functions.

\end{abstract}

\maketitle
\newfont{\bbmath}{msbm10 scaled\magstep1}
\newcommand{\C}{\mbox{\bbmath \symbol{'103}}}
\newcommand{\N}{\mbox{\bbmath \symbol{'116}}}
\newcommand{\Q}{\mbox{\bbmath \symbol{'121}}}
\newcommand{\R}{\mbox{\bbmath \symbol{'122}}}

\section{Introduction}

   It is well-known  \cite{titch} that the Riemann Hypothesis is equivalent to the 
  statement that all of
  the zeros of the  Riemann $\Xi(z)$ function are real.  
  $\Xi(z)$ is an even  entire
  function of $z$ with the integral representations
  \begin{equation}\label{xi}
	\Xi(z) = \int_{-\infty}^{\infty} \Phi(u)\,  e^{izu} \; du
	=2 \int_0^{\infty} \Phi(u) \cos(zu) \; du,
\end{equation}
where
 \begin{equation}\label{Phi}
	\Phi(u)= \sum^\infty_{n=1} (4 n^4\pi^2 e^{\frac{9}{2} u} 
	-  6 n^2\pi e^{\frac{5}{2} u}) e^{-n^2\pi e^{2u}}.
\end{equation}

In a 1926 paper  P\'olya  \cite{pol26}  observed  that 
   \begin{equation}\label{Phi-sim}
	\Phi(u) \sim 8 \pi^2  \cosh(\frac{9}{2} u)\,  e^{-2\pi \cosh(2u)}\quad \text{as} 
	\ u \rightarrow \pm \infty
\end{equation}
and, in view of this asymptotic equivalence  to $\Phi(u),$ 
considered the problem of determining
 whether or not  the entire function
\begin{equation}\label{xistar}
	\Xi^*(z) =  16 \pi^2 \int_0^{\infty} \cosh(\frac{9}{2} u)\,  e^{-2\pi \cosh(2u)} \cos(zu) \; du
\end{equation}
has only real zeros. Here, as is now  customary,  the capital letter $\Xi$ is used instead of  
the original lower case $\xi.$  P\'olya  was able to prove  that $\Xi^*(z)$ has only real zeros
by  using (in a  different notation)
a difference equation in  $z$ for the modified Bessel function of the third kind
 \cite{erd},   \cite{wat}
 \begin{equation}\label{K}
	K_{z}(a) = \int_0^{\infty} e^{-a \cosh u}\cosh(zu) \; du,\qquad a>0,
\end{equation}
to prove for each $ a>0$ that $K_{iz}(a)$   has only real zeros,
  and then applying the identity
  \begin{equation}\label{xiK}
	\Xi^*(z) = 4 \pi^2 [K_{\frac{1}{2}iz-\frac{9}{4}}(2 \pi)+
	K_{\frac{1}{2}iz+\frac{9}{4}}(2 \pi)]
	\end{equation}
and the special case $G(z)=K_{iz/2}(2 \pi), c=\frac{9}{4},$ of the
 lemma (derived via an infinite product representation for $G(z)$):
\begin{quote} {\bf Lemma.}
\em If $-\infty< c<\infty$
and $G(z)$ is an entire function of genus 0 or 1 that assumes real values
for real $z,$  has only real zeros and has at least one real zero,
then the function
$$G(z-ic)+G(z+ic)$$
also has only real zeros.
\end{quote}

More generally,
in a subsequent paper  P\'olya \cite{pol27} pointed out that  from the  case
$G(z)=K_{iz/2}(a)$ of this lemma it follows that each of the  entire functions
\begin{multline}\label{eq:F}
\quad \quad\ \,	F_{a,c}(z) = K_{i(z-ic)}(a)+ K_{i(z+ic)}(a)\\
=2 \int_0^{\infty}\cosh (cu)\,  e^{-a \cosh u}\cos(zu) \; du,\qquad a>0, \, -\infty< c<\infty,
	\end{multline}
 has only real zeros.
He also used a differential equation in the variable
$a$ to
give a proof of  the  reality of the
zeros of   $K_{iz}(a),  a>0,$  that was simpler than his previous proof.

    Our main aim in this paper is to show how  integrals of squares of 
 certain real-valued special functions can  be used to give new proofs of the reality
 of the zeros of the above  $\Xi^*(z), K_{iz}(a),$ and $F_{a,c}(z)$ functions.
 This paper is a sequel to  the author's 1994 paper  \cite{gg94}
in which he showed how sums of squares of certain real-valued special functions 
 could be used to prove the reality
 of the zeros of the Bessel functions $J_\alpha (z)$
when $\alpha \ge -1,$ confluent hypergeometric functions ${}_0F_1(c\/; z)$
when $c>0$ or $0>c>-1$, Laguerre polynomials $L_n^\alpha(z)$ when 
$\alpha \ge -2,$ Jacobi polynomials $P_n^{(\alpha,\beta)}(z)$
when $\alpha \ge -1$ and $ \beta \ge -1,$ and some other entire functions.
 Also see the applications of squares of real-valued special functions  in 
 \cite{ag76}, \cite{ag86}, \cite{de},  \cite{gg75},  \cite{gg75a}, \cite{gg77},  
 \cite{gg86},  \cite{gg89},    \cite{gg89a},   and \cite{gr}.

\section{ Reality of the zeros of the functions   $K_{iz}(a)$  when  $ a>0$}
\setcounter{equation}{0}

Let $a>0$ and  $z = x + i y,$ where
$x$ and $y$  are real variables. First observe that,
by the Meijer $G$-function representation for the
product of two  modified Bessel functions of the third kind
 in \cite[Eq. 5.6(66)]{erd}
and the definition of a Meijer $G$-function as a Mellin-Barnes integral  in
 \cite[Eq. 5.3(1)]{erd},
 \begin{multline}\label{K1}
\!\!	|K_{iz}(a)|^2 = K_{iz}(a)K_{i\bar{z}}(a) =
\frac{\sqrt\pi}{2}\, G_{24}^{40}\!\left[a^2\Big| \begin{matrix}0,\frac{1}{2} \\
ix,-ix,y,-y\end{matrix}\right] \  \\ 	
	=\frac{\sqrt\pi}{4\pi i}\int_{c-i\infty}^{c+i\infty}
	\frac{\Gamma(ix-s)\Gamma(-ix-s)\Gamma(y-s)\Gamma(-y-s)}
	{\Gamma(-s)\Gamma(\frac{1}{2}-s)} a^{2s}\, ds
\end{multline}
with $c<- |y|,$  where the path of integration is along the upwardly oriented
vertical line $\Re(s) = c.$
Next note that, by Gauss' summation formula  \cite[(1.2.11)]{gr},
$$	\sum_{k= 0}^\infty \frac{(y)_k (y)_k}{k!\,(y-s)_k} = 
	{}_2F_{1}\!\left[\begin{matrix}y,y\\
y-s \end{matrix};1 \right]= \frac{\Gamma(y-s)\Gamma(-y-s)}
	{\Gamma(-s)\Gamma(-s)} , \qquad \Re(s)<-y,$$
where $(y)_k =  \Gamma(y+k)/ \Gamma(y)$ and  the  series is absolutely convergent.
 Hence,
 $$	|K_{iz}(a)|^2  =
	\frac{\sqrt\pi}{4\pi i}\int_{c-i\infty}^{c+i\infty} \sum_{k= 0}^\infty
	\frac{((y)_k )^2\Gamma(ix-s)\Gamma(-ix-s)\Gamma(y-s)\Gamma(-s)}
	{k! \, \Gamma(\frac{1}{2}-s)\Gamma(y+k-s)} a^{2s}\, ds \quad$$
$$	= [K_{ix}(a)]^2+\frac{\sqrt\pi}{4\pi i}\sum_{k= 1}^\infty
	\frac{((y)_k )^2}{k!} 
	 \int_{c-i\infty}^{c+i\infty}
	\frac{\Gamma(ix-s)\Gamma(-ix-s)\Gamma(y-s)\Gamma(-s)}
	{\Gamma(\frac{1}{2}-s)\Gamma(y+k-s)} a^{2s} ds $$
and, since $(y)_k =y(y+1)_{k-1} $ for $k\ge 0,$
\begin{equation}
	|K_{iz}(a)|^2  = [K_{ix}(a)]^2 +y^2 L_a (x,y) \label{eq:Ky}
	\end{equation}
with
\begin{equation}\label{eq:L}
	 L_a (x,y)=\frac{\sqrt\pi}{4\pi i} \sum_{k= 0}^\infty
	\frac{((y+1)_k )^2}{(k+1)!} \!
	 \int_{c-i\infty}^{c+i\infty}
	\frac{\Gamma(ix-s)\Gamma(-ix-s)\Gamma(y-s)\Gamma(-s)}
	{\Gamma(\frac{1}{2}-s)\Gamma(y+k-s)} a^{2s} ds.
	\end{equation}
	
From \eqref{eq:Ky} it follows that in order to prove that $K_{iz}(a)$ 
	has only real zeros it suffices to prove that 
	\begin{equation}
	L_a (x,y) >0, \qquad  -\infty<x, y<\infty.  \label{eq:La}
	\end{equation} 
	  To prove	\eqref{eq:La} from \eqref{eq:L}
	  we observe that
	$$
	\int_0^1 t^{y-s-1} (1-t)^{k-1}\, dt = 
	\frac{\Gamma(k)\Gamma(y-s)}
	{\Gamma(y+k-s)}, \qquad k>0, \  \Re(s)<y,  $$
	   by the beta integral \cite[(1.11.8)]{gr}, and thus
	  	$$
	\frac{1}{2\pi i} \int_{c-i\infty}^{c+i\infty}
	\frac{\Gamma(ix-s)\Gamma(-ix-s)\Gamma(y-s)\Gamma(-s)}
	{\Gamma(\frac{1}{2}-s)\Gamma(y+k-s)} a^{2s} \, ds\qquad\qquad\qquad\qquad\quad\qquad $$
	 $$
=\frac{1}{\Gamma(k)} \int_0^1 t^{y-1} (1-t)^{k-1} \left[\int_{c-i\infty}^{c+i\infty}
\frac{\Gamma(ix-s)\Gamma(-ix-s)\Gamma(-s)}
	{2\pi i\, \Gamma(\frac{1}{2}-s)} \left(\frac{a^2}{t}\right)^s ds \right] dt $$
	$$ =\frac{2}{\sqrt{\pi}\Gamma(k)} \int_0^1 t^{y-1} (1-t)^{k-1}\,
	\left[K_{ix}\left(\frac{a}{\sqrt{t}}\right)\right]^2  dt \qquad\qquad\qquad\qquad\quad\qquad
\qquad$$
	  by Fubini's Theorem and (\ref{K1}), which gives
	   \begin{equation}\label{eq:Lt}
	 L_a (x,y) =\sum_{k= 0}^\infty
	\frac{((y+1)_k )^2}{k!(k+1)!} \int_0^1 t^{y-1} (1-t)^{k-1}\,
	\left[K_{ix}\left(\frac{a}{\sqrt{t}}\right)\right]^2  dt >0
\end{equation}
for 	 $ -\infty<x, y<\infty. $  Equations \eqref{eq:Ky} and \eqref{eq:Lt} can be combined
to give the  formula
    \begin{equation}\label{eq:Kp}
	 |K_{iz}(a)|^2  = [K_{ix}(a)]^2 + y^2 
	 \int_0^1 t^{y-1} {}_2F_{1}\!\left[\begin{matrix}y+1,y+1\\
2 \end{matrix};1-t \right]
	\left[K_{ix}\left(\frac{a}{\sqrt{t}}\right)\right]^2  dt, 
\end{equation}
from which it follows that $ K_{iz}(a)$ has only real zeros when $a>0,$
since the integrand is clearly nonnegative.

The reality of the zeros of $K_{iz}(a), a>0,$ can also be proved
by taking  the first and
 second partial derivatives  of the formula 
 \eqref{eq:Kp} with respect to $y$ to obtain the formulas
 \begin{equation}\label{eq:1st}
	y{\partial \over \partial y } |K_{iz}(a)|^2  = 
	 \int_0^1 y{\partial \over \partial y }f_t(y)\
	\left[K_{ix}\left(\frac{a}{\sqrt{t}}\right)\right]^2 \frac{ dt }{t}
\end{equation}
and
  \begin{equation}\label{eq:2nd}
	{\partial^2 \over \partial y^2} |K_{iz}(a)|^2  = 
	 \int_0^1{\partial^2 \over \partial y^2}f_t(y)\
	\left[K_{ix}\left(\frac{a}{\sqrt{t}}\right)\right]^2 \frac{ dt }{t}
\end{equation}
with
\begin{equation}\label{eq:f}
 f_t(y)= y^2 t^{y} {}_2F_{1}\!\left[\begin{matrix}y+1,y+1\\
2 \end{matrix};1-t \right] .
\end{equation}
In  \cite{ag97}  Askey and the author utilized the reality of the
zeros of the continuous dual Hahn polynomials \cite[p. 331]{aar} to show that
$f_t(y)$ (and a generalization of it) is an (even) absolutely monotonic function 
(one whose power series coefficients are nonnegative) of $y$ when $0<t<1,$ 
which shows that
\begin{equation}\label{eq:1pst}
y{\partial \over \partial y }f_t(y)\ge 0, \qquad -\infty<y<\infty, \  \  0<t<1, 
\end{equation}
and
\begin{equation}\label{eq:2pst}
{\partial^2 \over \partial y^2}f_t(y)\ge 0, \qquad -\infty<y<\infty, \  \  0<t<1.
\end{equation}
and, hence, that the integrals in 
 \eqref{eq:1st} and \eqref{eq:2nd} are positive when $y \ne 0.$   Thus,
it follows from  \eqref{eq:1st} that 
$|K_{iz}(a)|^2$  is an increasing (decreasing) even
function of $y$ when $y>0$ ($y<0),$ and it follows from \eqref{eq:2nd} that 
$|K_{iz}(a)|^2$
is a convex even non-constant function of $y,$ each of which implies that  for any $x$ 
the function $|K_{iz}(a)|^2$ assumes its minimum value
when $y=0$  and proves that $K_{iz}(a)$ has only real zeros.

Corresponding to the above proofs via formulas
\eqref{eq:1st} and \eqref{eq:2nd}, it should be noted that in 1913 Jensen \cite{jen}
 showed that each of the inequalities
\begin{equation}\label{Fx}
	\quad y{\partial \over \partial y }|F(x+i y)|^2 \ge 0, \quad -\infty<x, y<\infty,
\end{equation}
and
\begin{equation}\label{Fxy}
	{\partial^2 \over \partial y^2}|F(x+i y)|^2
	\ge 0,\quad -\infty<x, y<\infty,
\end{equation}
is necessary and sufficient for a real entire function $F(z) \not\equiv 0$ of 
genus 0 or 1  to have only real zeros
(see \cite[Chapter~2]{bo} and the  necessary and sufficient conditions in \cite{pol27a}).

\section{ Reality of the zeros of the functions  $\Xi^*(z)$ and $F_{a,c}(z)$ }
\setcounter{equation}{0}

Because of $\Xi^*(z)=4\pi^2 F_{2\pi,9/4}(z/2),$ it suffices to show how formula  \eqref{eq:Kp}
can be used to prove that $F_{a,c}(z)$ has only real zeros when $a, c>0.$ Fix $a, c>0$ 
and suppose that $z_0=x_0+i y_0$ is a zero of $F_{a,c}(z).$ Then 
$K_{i(x_0+i( y_0+c))}(a)=- K_{i(x_0+i( y_0-c))}(a)$ by \eqref{eq:F} and hence
\begin{multline}\label{eq:K0}
\qquad\quad 0=|K_{i(x_0+i( y_0+c))}(a) |^2-|K_{i(x_0+i( y_0-c))}(a) |^2 \\
= \int_0^1 [f_t(y_0+c)-f_t(y_0-c)]
	\left[K_{ix_0}\left(\frac{a}{\sqrt{t}}\right)\right]^2 \frac{ dt }{t} \qquad\qquad\qquad
\end{multline}
by \eqref{eq:Kp} and  \eqref{eq:f}.  Since $f_t(y)$ is an even convex non-constant function of $y$
when $0<t<1,$
\[f_t(y_0+c)-f_t(y_0-c) \left\{ \begin{array}{ll}
                                        >0 & \mbox{ if $ y_0>0, $}\\ 
                                        =0 & \mbox{ if $ y_0=0, $}\\ 
                                        <0 & \mbox{ if $ y_0<0 ,$}
                                        \end{array}
                                        \right. \]
 and  it follows from \eqref{eq:K0} that
$y_0=0$ and,  thus, the function
$F_{a,c}(z)$ has only real zeros when $a, c>0.$

One can also try to give another proof of the reality of the zeros of $F_{a,c}(z)$ 
via  a formula  for the function
\begin{multline}\label{eq:F2}
|F_{a,c}(z)|^2 =F_{a,c}(z)F_{a,c}(\bar{z})  \\
=K_{ix- y-c}(a) K_{ix+y-c}(a)+   K_{ix- y-c}(a)K_{ix+ y+c}(a) \qquad \\
+K_{ix- y+c}(a) K_{ix+y-c}(a)+   K_{ix- y+c}(a)K_{ix+ y+c}(a) \qquad\quad
\end{multline}
that contains integrals of nonnegative functions
as in   \eqref{eq:Kp}. By using the right-hand side of
\eqref{eq:F2},
\cite[Eq. 5.6(66)]{erd}, \cite[Eq. 5.3(1)]{erd}, 
Gauss' summation formula  and the beta integral, it can be shown that
\begin{multline}\label{eq:F3}
|F_{a,c}(z)|^2 =[F_{a,c}(x)]^2  
+\int_0^ 1 f_t(y) [F_{a/\sqrt{t}, \, c}(x)]^2 \frac{ dt }{t}  
+ \int_0^1 g_{t, c}(y) |K_{i(x+ic)}(a/\sqrt{t})|^2 \frac{ dt }{t} 
\end{multline}
with
\begin{multline}\label{eq:g}
g_{t, c}(y)=y(y+2 c) t^{y+c} {}_2F_{1}\!\left[\begin{matrix}y+1,y+2 c+1\\
2 \end{matrix};1-t \right] \\
+y(y-2 c) t^{y-c} {}_2F_{1}\!\left[\begin{matrix}y+1,y-2 c+1\\
2 \end{matrix};1-t \right] 
-2 y^2 t^{y} {}_2F_{1}\!\left[\begin{matrix}y+1,y+1\\
2 \end{matrix};1-t \right] .
\end{multline}
The function $g_{t, c}(y)$ is clearly an even function of $c$ such that
$g_{t, 0}(y)=0$ and $g_{t, c}(0)=0.$
It can be shown that $g_{t, c}(y)$ is also an even function of $y$ by applying
the Euler transformation formula  \cite[(1.4.2)]{gr}  to the hypergeometric
functions in \eqref{eq:g}.
Extensive analysis of $g_{t, c}(y)$ via Mathematica  and generating
functions of the form  \cite[(5)]{ag97}  strongly suggest that
$g_{t, c}(y)$ is   nonnegative, convex,  and an absolutely monotonic function 
of $y$ when $0<t<1$  and $ -\infty<c<\infty.$  If the nonnegativity, convexity, 
or absolute monotonicity of $g_{t, c}(y)$ could be proved, then 
\eqref{eq:F3} and its partial derivatives with respect to 
 $y$ would give additional proofs of the 
reality of the zeros of $F_{a,c}(z).$

 P\'olya \cite{pol27} derived a theorem concerning the zeros of  Fourier transforms
  and universal factors, and then used it to prove that his  \cite[p. 317]{pol26}
 \begin{equation}\label{eq:xistarstar}
	\Xi^{**}(z) = 8 \pi
	 \int_0^{\infty} \left[ 2 \pi \cosh(\frac{9}{2} u)-3\cosh(\frac{5}{2} u)\right]
	  e^{-2\pi \cosh(2u)} \cos(zu) \; du,
\end{equation}
function and the more general functions
\begin{equation}\label{eq:xiA}
	\Xi_{A,B,a,b,c}(z) = 
	 \int_0^{\infty} [ A \cosh(a u)-B\cosh(b u)]
	  e^{-c \cosh(u)} \cos(zu) \; du
\end{equation}
with $A>B>0,  a>b>0,  c>0,$ have only real zeros.  Analogous
to  the formulas in \eqref{eq:Kp} and \eqref{eq:F3},
it might be possible to derive  formulas containing 
squares of real-valued functions
that give new proofs of the reality of the 
zeros of the  entire functions in \eqref{eq:xistarstar}
and \eqref{eq:xiA}, and even, perhaps, prove that some of
the  entire functions in
Hejhal \cite[(0.3)]{h} have only  real zeros.

 \section*{Acknowledgment}  
 The author wishes to thank the referee for suggesting some improvements
 in the paper.

\end{document}